\documentstyle[amssymb,amssymb,draft,12pt]{amsart}

\newcommand{\wt}{\widetilde}
\newcommand{\dnm}{{\delta}_{n, -m}}
\newcommand{\Z}{\Bbb Z}
\newcommand{\La}{\Lambda}

\newcommand{\son}{t - t^{-1}}
\newcommand{\sonn}{t^{-1} - t}
\newcommand{\stw}{t^{2} - t^{-2}}
\newcommand{\sth}{t^{3} - t^{-3}}
\newcommand{\sfo}{t^{4} - t^{-4}}
\newcommand{\sfi}{t^{5} - t^{-5}}
\newcommand{\ssi}{t^{6} - t^{-6}}
\newcommand{\sei}{t^{8} - t^{-8}}
\newcommand{\con}{t + t^{-1}}
\newcommand{\ctw}{t^{2} + t^{-2}}
\newcommand{\cth}{t^{3} + t^{-3}}
\newcommand{\csi}{t^{6} + t^{-6}}
\newcommand{\den}{t^{n-1} + t^{-(n-1)}}
\newcommand{\si}{t^i - t^{-i}}
\newcommand{\g}{{\frak g}}
\newcommand{\hg}{\hat{{\frak g}}}
\newcommand{\UU}{U_q(\hat{\frak g})}
\newcommand{\UUU}{\tilde{U}_q(\hat{\frak g})}
\newcommand{\ZZ}{Z_q(\hat{\frak g})}
\newcommand{\HH}{H_q({\frak g})}
\newcommand{\WW}{W_q({\frak g})}

\newtheorem{th}{Theorem}[section]

\newtheorem{prop}[th]{Proposition}

\title[Deformations of $W$-algebras]
{Deformations of the classical ${\cal W}$-algebras 
associated to $D_n, E_6$ and $G_2$}
\date\today
	\address{Department of Mathematics, University of California, Berkeley, 94720}
\email{kogan@@math.berkeley.edu}

\author{Alexander Kogan}
\begin{document}

\maketitle

\section{introduction.}
The purpose of this paper is the computation of the Poisson brackets in
the deformed $W$-algebras $\WW$, where $\g$ is of type $D_n, E_6$ or
$G_2$. Let us first briefly recall some facts about two main descriptions of 
$W$-algebras. 

The first description is via the Drinfeld-Sokolov reduction.
Let ${\frak g}$ be a simple Lie algebra, and ${\frak n}$ be its nilpotent
subalgebra. Let $M = \mu^{-1}(f)/LN$, where $\mu : \hg^* \rightarrow
L{\frak n}^*$ is the momentum map and $f \in L{\frak n}^*$ is a certain
character \cite{Dr-So}. Then $M$ is a Poisson manifold, and the ordinary
$W$-algebra $W(\g)$ is the Poisson algebra of functions on $M$.
The manifold $M$ can be identified with the space of certain differential
operators. For example if $\g = {\frak sl}_n$ then these operators are of the
form
 \[ \partial^n + a_{n-2}\partial^{n-2}+ \cdots + a_0,\]
 and the Poisson structure under consideration is called the
second Gelfand-Dickey bracket \cite{Ge-Di}. Recently in 
\cite{Fr-Re-Se,Se-Se} a $q$-deformation of $W$-algebras was obtained by 
considering the space of $q$-difference operators. 

Let us now recall the second description of $W$-algebras. It was proved by 
B.~Feigin and E.~Frenkel \cite{Fe-Fr,Fr} that as Poisson algebras 
$W(\g^L)$ is isomorphic to the center $Z(\hat{\frak g})$
 of the completion of the universal enveloping algebra
 $\tilde{U}(\hat{\frak g})_{-h^\vee}$ of $\hg$ 
at the critical level $-h^\vee$, 
where $h^\vee$ is the dual Coxeter number and $\g^L$ is the Langlands
dual of the algebra $\g$. This description was used in 
\cite{Fr-Re:algebras} to produce a $q$-deformed $W$-algebras
$\WW$, where $\g$ is of the classical type, and to compute the 
Poisson structure in the $A_n$ case. In the cases $B_n$ and $C_n$ some
Poisson brackets were computed in \cite{Fr-Re:BC}. 

It is convenient to study $W$-algebras via the Miura transformation. For
ordinary $W$-algebras it can be defined as follows
(see \cite{Fe-Fr,Fr,Fr-Re:algebras}). One considers 
the Wakimoto homomorphism from $\tilde{U}(\hg)$ to the tensor 
product of some Heisenberg algebra and some commutative algebra $H(\g)$ - 
algebra of functionals on some hyperplane in the dual space to the 
Heisenberg subalgebra $\hat{\frak{h}}$ of $\hg$. The restriction of this map 
to the center gives the homomorphism $Z(\hg) \rightarrow H(\g)$, composition 
of which with the isomorphism $W(\g^L) \simeq Z(\hg)$ on the left is just the 
Miura transformation \cite{Dr-So,Fr-Re:algebras}. 

The $q$-deformed version of this picture is the Wakimoto realization of 
$\UUU$ in the tensor product of a certain Heisenberg algebra and some 
Heisenberg-Poisson algebra $\HH$.
The restriction to the center $\ZZ$ 
of $\UUU$ gives the $q$-deformed Miura 
transformation $\ZZ \rightarrow \HH$. The image is called the
 $q$-$W$-algebra $\WW$.
Thus, in order to describe deformed $W$-algebras we have to describe the
Heisenberg-Poisson algebra $\HH$ and the generators of $\WW$. In
\cite{Fr-Re:algebras} E.~Frenkel and \\ N.~Reshetikhin did this in the
$A_n$ case using the explicit formulas for the Wakimoto realization
\cite{AOS}. Motivated by these results they gave a conjectural
description of $\HH$ for general ${\frak g}$ (see Sect. 11 of
\cite{Fr-Re:algebras} and the next section) and of $\WW$ for ${\frak g}$
of classical series. The key element of this conjecture was that the
formulas for the generators of the deformed $W$-algebra coincide with the
formulas for the eigenvalues of the corresponding transfer-matrices 
obtained by  analytic Bethe Ansatz (see Conjecture 1 of
\cite{Fr-Re:algebras}). In order to verify this conjecture, one has to
check that the Poisson brackets between the generators of $\WW$,
constructed in this way, close among themselves. This had been done in
\cite{Fr-Re:algebras} for the $A_n$ series and in \cite{Fr-Re:BC,
 Fr-Re:def} for the
$B_n$ and $C_n$ series. However the question remained open for other
series.

In this paper we study the case when $\g$ is of the type $D_n, E_6$ or
$G_2$. We exhibit the generators and relations of the algebra $\HH$
explicitely and compute the Poisson brackets between them. Next, we
construct the generators of the $W$-algebra $\WW$, following the conjecture
of \cite{Fr-Re:algebras} that they should coincide with the
corresponding formulas for the eigenvalues of transfer-matrices
(see \cite{Resh}). Finally, we compute the Poisson brackets between them.

The paper is organized as follows. In section 2 we describe the algebra $\HH$
via the generators and relations. Sections 3, 4 and 5 are devoted to the cases of $D_n, E_6$ and $G_2$ respectively. Each of these sections is divided in three subsections. In the first one we define matrices which are used in the 
construction of $\HH$. In the second subsection we describe the new set of 
generators of $\HH$ which is convenient for the computation of the Poisson
brackets. Then we axiomatically define the
 generators of the corresponding
$W$-algebra. In the last subsection we compute the Poisson brackets between 
the generators of $\HH$ and $\WW$. 
 
{\bf Acknowledgements.} I am very grateful to E.~Frenkel and \\ 
N.~Reshetikhin for useful ideas and discussions. I would like to thank RIMS, and particularly Tetsuji Miwa,
for the hospitality in July-August of 1997, when this work was 
being completed.

\section{Heisenberg-Poisson algebras.}

In this section we describe the algebra $\HH$. The matrices $M, \tilde{M}$ and
$D$ will be explicitely given in the next sections. It should be noted that 
$\tilde{M}$ is the deformation of twice the 
symmetrized Cartan matrix of the corresponding 
algebra. The presentation in this section follows \cite{Fr-Re:BC,
Fr-Re:def}.

Let $\g$ be either $D_n, E_6$ or $G_2$. Let $\UUU$ be the completion of the quantum  universal enveloping algebra $\UU$ of $\g$ as defined in
 \cite{Fr-Re:algebras}. 

We consider a Heisenberg-Poisson algebra $\HH$ with generators $a_i[n], \\
 1 \leq i \leq $ rank($\g$) and 
relations: \[ \{a_i[n], a_j[m]\} = \wt{M}_{ij}(q^n)\dnm. \]
There is unique set of "dual" generators $y_i[n]$ such that
\[ \{y_i[n], a_j[m]\} = D_{ij}(q^n)\dnm. \]
Then $y_i[n]$ satisfy 
\[ \{y_i[n], y_j[m]\} = M_{ij}(q^n)\dnm. \]
Let's form the generating series:
\[ Y_i(z) = q^{-2(\rho, \omega_i)} \exp \left( -\sum_{m \in \Z} y_i[m]z^{-m} \right) \]
They satisfy the following relations:
\[ \{Y_i(z), Y_j(w)\} = {\cal M}_{ij} \left( \frac{w}{z} \right) Y_i(z)Y_j(w), \]
where 
\[ {\cal M}_{ij}(x) = \sum_{m \in \Z} M_{ij}(q^m)x^{m}.\]

The coefficients of the generationg functions $Y_i(z)$ generate the algebra
$\HH$. In the next sections we will introduce the new generating functions 
$\La_i(z)$ which also have generators as coefficients. 
Finally the generating functions
whose coefficients generate $\WW$ will be denoted $T_i(z)$, where 
$1 \leq i \leq$ (rank of the Cartan subalgebra of $\g$). In the $D_n$ 
\cite{Fr-Re:BC} and $G_2$ cases all $T_i(z)$'s can be constructed explicitely,
 whereas in the the case of $E_6$ we explicitely construct only $T_1(z)$. 

\section{$D_n$ case.}
\subsection{Matrices}

Consider the matrices $M(t), D(t), \wt{M}(t)$ defined as follows. Let $M(t) = (M_{ij}(t)), 1 \leq i,j \leq n$, where
\begin{equation*}
\begin{array}{l}
\displaystyle M_{ij}(t) = \displaystyle \frac{(t^{min(i,j)} - t^{-min(i,j)})(t^{n - 1 - max(i,j)} + t^{-(n - 1 - max(i,j)})}{(\den)}, \\
\displaystyle M_{ni}(t) = M_{n-1,i}(t) = \displaystyle \frac{(\si)}{(\den)}, 1 \leq i,j \leq n-2, \\
\displaystyle M_{n,n-1}(t) = \displaystyle \frac{(t^{n-2} - t^{-(n - 2})}{(\con)(\den)}, \\
\displaystyle M_{n-1,n-1}(t) = M_{nn}(t) = \displaystyle \frac{(t^{n} - t^{-n})}{(\con)(\den)}. 
\end{array}
\end{equation*}

Let $D(t) = (t - t^{-1})\cdot I_{n}$, where $I_{n}$ is the $n \times n$ identity matrix. Then  $$\wt{M}(t) = D(t)M(t)^{-1}D(t) = $$
$$\begin{pmatrix}
\stw & -(\son) & \hdots & 0 & 0 & 0 \\
-(\son) & \stw & \hdots & 0 & 0 & 0 \\
\hdotsfor{6} \\
0 & 0 & \hdots & \stw & -(\son) & -(\son) \\
0 & 0 & \hdots & -(\son) & \stw & 0 \\
0 & 0 & \hdots & -(\son) & 0 & \stw
\end{pmatrix},$$ which is a $t$--deformation of twice the Cartan
matrix of $D_n$.

\subsection{Generators}
Introduce as in \cite{Fr-Re:algebras}
\begin{align*}
\La_i(z) = & Y_i(zq^{-i+1}) Y_{i-1}^{-1}(zq^{-i}), \quad \quad
i=1,\ldots,n-2, \\
\La_{n-1}(z) = & Y_n(zq^{-n+2}) Y_{n-1}(zq^{-n+2})
Y_{n-2}^{-1}(zq^{-n+1}), \\
\La_{n}(z) = & Y_{n-1}(zq^{-n+2}) Y_n^{-1}(zq^{-n}), \\
\La_{n+1}(z) = & Y_{n}(zq^{-n+2}) Y_{n-1}^{-1}(zq^{-n}), \\
\La_{n+2}(z) = & Y_{n-2}(zq^{-n+1}) Y_{n-1}^{-1}(zq^{-n}) Y_n^{-1}(zq^{-n}), \\
\La_{2n-i+1}(z) = & Y_{i-1}(zq^{-2n+i+2}) Y_i^{-1}(zq^{-2n+i+1}) \quad
\quad i=1,\ldots,n-2,
\end{align*}
where $Y_0(z) = 1$.

{\em Remark.}  The relations between $Y_i(z)$ and the functions 
$Q_i(u), 1 \leq i \leq n,$ which appear in \cite{Resh} in
the formulas for the eigenvalues of the transfer matrices of the 
$D^{(1)}_n$ model are as follows:
\[ Y_i(zq^m) = \frac{Q_i(u+\frac{m+1}{2}\eta)}{Q_i(u+\frac{m-1}{2}\eta)}\]
Let
\begin{align*}
 T_1(z) = & \sum_{i=1}^{2n} \La_i(z).\\
 T_2(z) = & \sum_{(i,j) \in S} \La_i(z)\La_j(zq^2),  
\end{align*}
    where the set $S$ consists of pairs $(i,j)$ such that either
    $i<j$ or $(i,j) = (n+1,n).$ 
\begin{align*}
\{ \La_i(z),\La_i(w) \} &= {\cal M}_{11}\left( \frac{w}{z} \right) \La_i(z)
\La_i(w), \\
\{ \La_i(z),\La_j(w) \} &= {\cal M}_{11}\left( \frac{w}{z} \right) \La_i(z)
\La_j(w) +
\left( \delta \left( \frac{w}{zq^2} \right) - \delta \left( \frac{w}{z}
\right) \right) \\ &+ \left( \delta \left( \frac{w}{zq^{2n-2i}} \right)
\delta_{i+j,2n+1}-\delta \left( \frac{w}{zq^{2n-2i-2}} \right)
\delta_{i+j,2n+1} \right) \\
&\times \La_i(z) \La_j(w),
\end{align*}
if $i<j$. 
\begin{align*}
\{ T_1(z),T_1(w) \} &= {\cal M}_{11} \left( \frac{w}{z} \right) T_1(z)
T_1(w) \\ &+ \delta \left( \frac{w}{zq^2} \right) T_2(z) - \delta \left(
\frac{wq^2}{z} \right) T_2(w) \\ &+ \delta \left( \frac{w}{zq^{2n-2}}
\right) - \delta \left( \frac{w^{2n-2}}{z} \right).
\end{align*}
\section{$E_6$ case.}
\subsection{Matrices}
Consider the matrices $M(t), D(t), \wt{M}(t)$ defined as follows. Let $M(t) = (M_{ij}(t)), 1 \leq i,j \leq 6$, $M_{ij}(t) = M_{ji}(t)$, where
\begin{align*}
\displaystyle M_{11}(t) = M_{55}(t) = & \displaystyle \frac{(\son)(\sei)}{(\csi)(\sth)} \\ 
\displaystyle M_{12}(t) = M_{45}(t) = & \displaystyle \frac{(\son)(\sfi)(\ctw)}{(\csi)(\sth)} \\   
\displaystyle M_{22}(t) = M_{44}(t) = & \displaystyle \frac{(\sfo)(\sfi)}{(\csi)(\sth)} \\
\displaystyle M_{13}(t) = M_{26}(t) = M_{46} = M_{35} = & \displaystyle \frac{(\sfo)}{(\csi)} \\ 
\displaystyle M_{23}(t) = M_{34}(t) = & \displaystyle \frac{(\sfo)(\con)}{(\csi)} \\
\displaystyle M_{33}(t) = & \displaystyle \frac{(\sth)(\con)(\ctw)}{(\csi)} \\
\displaystyle M_{16}(t) = M_{56}(t) = & \displaystyle \frac{(\son)(\ctw)}{(\csi)} \\
\displaystyle M_{36}(t) = & \displaystyle \frac{(\sth)(\ctw)}{(\csi)} 
\end{align*}

\begin{align*}
\displaystyle M_{66}(t) = & \displaystyle \frac{(\sfo)(\cth)}{(\con)(\csi)} \\
\displaystyle M_{14}(t) = M_{25}(t) = & \displaystyle \frac{(\stw)(\sfo)}{(\csi)(\sth)} \\
\displaystyle M_{24}(t) = & \displaystyle \frac{(\stw)(\sfo)(\con)}{(\csi)(\sth)} \\
\displaystyle M_{15}(t) = & \displaystyle \frac{(\son)(\sfo)}{(\csi)(\sth)}. 
\end{align*}
Let $D(t) = (t - t^{-1})\cdot I_{6}$, where $I_{6}$ is the $6 \times 6$ identity matrix. Then  $$\wt{M}(t) = D(t)M(t)^{-1}D(t) = $$
\[ \left( \begin{array}{cccccc}
\stw & \sonn & 0 & 0 & 0 & 0 \\
\sonn & \stw & \sonn & 0 & 0 & 0 \\
0 & \sonn & \stw & \sonn & 0 & \sonn \\
0 & 0 & \sonn & \stw & \sonn & 0 \\
0 & 0 & 0 & \sonn & \stw & 0 \\
0 & 0 & \sonn & 0 & 0 & \stw 
\end {array} \right) \]
is a $t$-deformation of twice the Cartan matrix of $E_6$.

\subsection{Generators}

Introduce
\begin{align*}
\La_{1}(z) = & Y_{1}^{-1}(zq^{-8}) Y_{2}(zq^{-7}) Y_{3}^{-1}(zq^{-8}) Y_{6}(zq^{-7}), \\
\La_{2}(z) = & Y_{1}^{-1}(zq^{-8}) Y_{2}(zq^{-7}) Y_{6}^{-1}(zq^{-9}), \\
\La_{3}(z) = & Y_{1}^{-1}(zq^{-8}) Y_{3}(zq^{-6}) Y_{4}^{-1}(zq^{-7}), \\
\La_{4}(z) = & Y_{1}^{-1}(zq^{-8}) Y_{4}(zq^{-5}) Y_{5}^{-1}(zq^{-6}), \\
\La_{5}(z) = & Y_{2}^{-1}(zq^{-9}) Y_{3}(zq^{-8}) Y_{6}^{-1}(zq^{-9}), \\
\La_{6}(z) = & Y_{2}^{-1}(zq^{-9}) Y_{6}(zq^{-7}), \\
\La_{7}(z) = & Y_{3}^{-1}(zq^{-10}) Y_{4}(zq^{-9}), \\
\La_{8}(z) = & Y_{4}^{-1}(zq^{-11}) Y_{5}(zq^{-10}), \\
\La_{9}(z) = & Y_{1}(zq^{-6}) Y_{2}^{-1}(zq^{-7}) Y_{3}(zq^{-6}) Y_{4}^{-1}(zq^{-7}),\\ 
\La_{10}(z) = & Y_{1}(zq^{-6}) Y_{2}^{-1}(zq^{-7}) Y_{4}(zq^{-5}) Y_{5}^{-1}(zq^{-6}), \\
\La_{11}(z) = & Y_{1}(zq^{-6}) Y_{3}^{-1}(zq^{-8}) Y_{6}(zq^{-7}), \\
\La_{12}(z) = & Y_{1}(zq^{-6}) Y_{6}^{-1}(zq^{-9}),\\ 
\La_{13}(z) = & Y_{2}(zq^{-5}) Y_{3}^{-1}(zq^{-6}) Y_{4}(zq^{-5}) Y_{5}^{-1}(zq^{-6}), \\
\La_{14}(z) = & Y_{2}(zq^{-5}) Y_{4}^{-1}(zq^{-7}),\\ 
\La_{15}(z) = & Y_{3}(zq^{-4}) Y_{5}^{-1}(zq^{-6}) Y_{6}^{-1}(zq^{-5}),\\
\La_{16}(z) = & Y_{5}^{-1}(zq^{-6}) Y_{6}(zq^{-3}), \\
\La_{17}(z) = & Y_{1}^{-1}(zq^{-8}) Y_{5}(zq^{-4}), \\
\La_{18}(z) = & Y_{1}(zq^{-6}) Y_{2}^{-1}(zq^{-7}) Y_{5}(zq^{-4}), \\
\La_{19}(z) = & Y_{2}(zq^{-5}) Y_{3}^{-1}(zq^{-6}) Y_{5}(zq^{-4}), \\
\La_{20}(z) = & Y_{3}(zq^{-4}) Y_{4}^{-1}(zq^{-5}) Y_{5}(zq^{-4}) Y_{6}^{-1}(zq^{-5}), \\
\La_{21}(z) = & Y_{4}^{-1}(zq^{-5}) Y_{5}(zq^{-4}) Y_{6}(zq^{-3}), \\
\La_{22}(z) = & Y_{4}(zq^{-3}) Y_{6}^{-1}(zq^{-5}), \\
\La_{23}(z) = & Y_{3}^{-1}(zq^{-4}) Y_{4}(zq^{-3}) Y_{6}(zq^{-3}), \\
\La_{24}(z) = & Y_{2}^{-1}(zq^{-3}) Y_{3}(zq^{-2}), \\
\La_{25}(z) = & Y_{1}^{-1}(zq^{-2}) Y_{2}(zq^{-1}), \\
\La_{26}(z) = & Y_{1}(z), \\
\La_{27}(z) = & Y_{5}^{-1}(zq^{-12}).
\end{align*}
{\em Remark.}  The relations between $Y_i(z)$
 and the functions $Q_i(u), 1\leq i \leq 6,$ which appear in \cite{Resh} in
the formulas for the eigenvalues of the transfer matrices of the 
$E^{(1)}_6$ model are as follows:
 \[Y_i(zq^m) = \frac{Q_{\sigma(i)}(u-(m+1)\eta)}{Q_{\sigma(i)}(u-(m-1)\eta)}, \]
where $\sigma$ is a permutation $(1)(2)(3)(456)$.

Let \[T_1(z) = \sum_{i=1}^{27} \La_i(z).\]

\subsection{Poisson brackets}

\begin{equation*}
\begin{array}{ll}
\{\La_i(z), \La_j(w)\} = & M_{11} \left( \frac{w}{z}
                                \right)\La_i(z)\La_j(w)\\ 
                         & + (\text{sum of }\delta \text{- functions})\La_i(z)\La_j(w).
\end{array}
\end{equation*}

\begin{align*}
\{ T_1(z),T_1(w) \} &=  M_{11} \left( \frac{w}{z} \right) T_1(z)
T_1(w) \\ &+ \delta \left( \frac{wq^2}{z} \right) T_2(z) - \delta \left(
\frac{w}{zq^2} \right) T_2(w) \\ &+ \delta \left( \frac{wq^8}{z}
\right) T_5(zq^4) - \delta \left( \frac{w}{zq^8} \right) T_5(wq^4),
\end{align*}
where $T_5(z)$ is the $W$-algebra generator corresponding to the fifth 
fundamental weight (which is dual to the first one in our notation).

\section{$G_2$ case.}
\subsection{Matrices}

Consider the matrices $M(t), D(t)$ and $\wt{M}(t)$ defined as follows. Let
$M(t) = (M_{ij}(t)), 1 \leq i,j \leq 2$, where
\begin{equation*} 
\begin{array}{l}
\displaystyle M_{22}(t) = \frac{(\sth)(\con)(\ctw)}{\csi} \\
\displaystyle M_{11}(t) = \frac{(\cth)(\son)(\ctw)}{\csi}\\
\displaystyle M_{12}(t) = M_{21}(t) = \frac{(\sth)(\ctw)}{\csi}.
\end{array}
\end{equation*}
Let \[D(t) = \left( \begin{array}{cc}
\son & 0 \\
0 & \sth
                    \end{array} \right) \]
Then $$\wt{M}(t) = D(t)M(t)^{-1}D(t) = $$
\[ \left( \begin{array}{cc}
\stw & -(\sth)\\
-(\sth) & \ssi
           \end{array} \right) \] 
is a $t$-deformation of the symmetrized Cartan matrix of $G_2$.

\subsection{Generators}
Introduce
\begin{equation*}
\begin{array}{l}
\La_1(z) = Y_1(z), \\
\La_2(z) = Y_1^{-1}(zq^{-2})Y_2(zq^{-1}), \\
\La_3(z) = Y_1(zq^{-4})Y_1(zq^{-6})Y_2^{-1}(zq^{-7}), \\
\La_4(z) = Y_1(zq^{-4})Y_1^{-1}(zq^{-8}), \\
\La_5(z) = Y_1^{-1}(zq^{-6})Y_1^{-1}(zq^{-8})Y_2(zq^{-5}), \\
\La_6(z) = Y_1(zq^{-10})Y_2^{-1}(zq^{-11}), \\
\La_7(z) = Y_1^{-1}(zq^{-12}).
\end{array} 
\end{equation*}
{\em Remark.}  The relations between $Y_i(z)$ and the functions 
$Q_i(u), i=1,2,$ which appear in \cite{Resh} in
the formulas for the eigenvalues of the transfer matrices of the 
$G^{(1)}_2$ model are as follows:
\[ Y_1(zq^m) = \frac{Q_1(u+\frac{13+m}{3}\eta)}{Q_1(u+\frac{11+m}{3}\eta)}\]
\[ Y_2(zq^m) = \frac{Q_2(u+\frac{15+m}{3}\eta)}{Q_2(u+\frac{9+m}{3}\eta)}\]
.

Let \[T_1(z) = \sum_{i=1}^7 \La_i(z).\]
    \[T_2(z) = \sum_{i=2}^7 \La_1(z)\La_i(zq^{2})
             + \sum_{i=2}^6 \La_i(z)\La_7(zq^{2})+\] 
	    \[+\La_2(z)\La_5(zq^{2})+
	    \La_2(z)\La_6(zq^{2})+\]
	    \[+\La_3(z)\La_5(zq^{2})+\La_3(z)\La_6(zq^{2}).\]

\subsection{Poisson brackets}

\begin{equation*}
\begin{array}{ll}
\{\La_i(z), \La_j(w)\} = & M_{11} \left( \frac{w}{z} \right)\La_i(z)\La_j(w)\\ 
                 & +(\text{sum of }\delta \text{-functions}) \La_i(z)\La_j(w).\\
\end{array}
\end{equation*}

\begin{align*}
\{ T_1(z),T_1(w) \} &= M_{11} \left( \frac{w}{z} \right) T_1(z)
T_1(w) \\ &+ \delta \left( \frac{w}{zq^2} \right) T_2(z) - \delta \left(
\frac{wq^2}{z} \right) T_2(w) \\ &+ \delta \left( \frac{w}{zq^8}
\right) T_1(zq^4) - \delta \left( \frac{wq^8}{z} \right) T_1(wq^4)\\
 &+ \delta \left( \frac{w}{zq^{12}} \right) - \delta \left( 
\frac{wq^{12}}{z} \right).
\end{align*}

There is an obvious
\begin{prop}
In the cases $E_6$ and $G_2$ 
if we replace $Y_{i}^{\epsilon}(zq^n) \text{ by } \\
 Y_{i}^{-\epsilon}(zq^{-n})$ in $T_1(z)$ then we obtain
$T_\alpha(zq^{12})$, where $T_\alpha(z)$ corresponds to the dual root (i.e. 
$\alpha = 5 \text{ in } E_6$ case and $1 \text{ in } G_2$ case).
\end{prop}

\end{document}